\documentclass[reqno]{amsart}
\usepackage{amsmath}

\def\eqn#1{(\ref{eq:#1})}
\newcommand{\qBin}[3]{\genfrac{[}{]}{0pt}{0}{#1}{#2}_{#3}}

\newcommand{\bin}{\binom}
\newcommand{\beqs}{\begin{equation*}}
\newcommand{\eeqs}{\end{equation*}}
\newcommand{\beq}{\begin{equation}}
\newcommand{\eeq}{\end{equation}}

\theoremstyle{plain}
\newtheorem{thm}{Theorem}

\numberwithin{equation}{section}
\allowdisplaybreaks[2]

\begin{document}

\title[\tiny{The tri-pentagonal number theorem and related identities}]
      {The tri-pentagonal number theorem and related identities}

\author[Alexander~Berkovich]{Alexander Berkovich}
\address{Department of Mathematics, 
         University of Florida,
         Gainesville, FL~32611}
\email{alexb@math.ufl.edu}
\thanks{Research was supported in part by NSA grant H98230-07-01-0011.}


	           
\begin{abstract}
I revisit an automated proof of Andrews' pentagonal number theorem found by Riese. 
I uncover a simple polynomial identity hidden behind his proof. 
I explain how to use this identity to prove Andrews' result along with a variety of new formulas of similar type. 
I reveal an interesting relation between the tri-pentagonal theorem and items (19), (20), (94), (98) 
on the celebrated Slater list. Finally, I establish a new infinite family of multiple series identities.
\end{abstract}

\maketitle
\section{\bf Introduction}
\label{sec:1}
\medskip

The Gaussian or $q$-binomial coefficients are polynomials in $q$ defined by 
\begin{align*}
\qBin{n+m}{n}{q}:=
\begin{cases} \frac{(q)_{n+m}}{(q)_n(q)_m}, & \mbox{if } n,m\in\mathbb{N}, \\
               0, & \mbox{otherwise}. 
\end{cases}
\end{align*}
Here $(q)_n=\prod^n_{j=1}(1-q^j)$. We shall require the more general $q$-shifted factorials defined by
\begin{align*}
(a)_n=(a;q)_n:=
\begin{cases}  1, & \mbox{ if } n=0, \\
               \prod_{j=0}^{n-1}(1-aq^j) & \mbox{ if } n>0, \\
               \prod_{j=1}^{-n}\frac{1}{1-aq^{-j}} & \mbox{ if } n<0. 
\end{cases}
\end{align*}
We note that
\beqs
\frac{1}{(q)_n}=0, \mbox{ if } n<0,
\eeqs
and that
\beqs
\lim_{L\rightarrow\infty}\qBin{L}{j}{q}=\frac{1}{(q)_j}.
\eeqs
Here and hereafter $|q|<1$.
We shall also use the following notations
\beqs
(a_1,a_2,\ldots,a_k)_n=(a_1,a_2,\ldots,a_k;q)_n:=\prod_{i=1}^k(a_i)_n,
\eeqs
\beqs
(a)_\infty=(a;q)_\infty:=\lim_{n\rightarrow\infty}(a)_n,    
\eeqs
\beqs 
[z;q]_\infty:=(z,\frac{q}{z})_\infty.
\eeqs
The literature on $q$-series abounds with numerous identities of the type
\beq
\sum_{i=-\infty}^\infty A(i,q)\qBin{2L}{L-i}{q}=B(L,q).
\label{eq:1.1}
\eeq
For example,
\begin{align}
\sum_{j=-\infty}^\infty (-1)^jq^{\bin{j}{2}}\qBin{2L}{L-j}{q} & =\delta_{L,0}, 
\label{eq:1.2} \\
\sum_{j=-\infty}^\infty (-1)^jq^{\bin{j}{2}}\qBin{2L}{L-j}{q^2} & =(-1)^Lq^{L^2}(q;q^2)_L, 
\label{eq:1.3} \\
\sum_{j=-\infty}^\infty (-1)^jq^{2\bin{j}{2}}\qBin{2L}{L-j}{q} & =q^L(q;q^2)_L.
\label{eq:1.4}
\end{align}
Here $\delta_{i,j}$ is the Kronecker delta and $\bin{i}{2}:=\frac{i(i-1)}{2}$. \\
The observant reader might have recognized \eqn{1.2} as a special case of
\beq
\sum_{j=-\infty}^\infty (-1)^j z^j q^{\bin{j+1}{2}}\qBin{2L+a}{L-j}{q}=(\frac{1}{z})_{L+a}(qz)_L,
\label{eq:1.5}
\eeq
where $a\in\mathbb N$. The above formula is due to Cauchy.
It is a finite form of the celebrated Jacobi triple product identity 
\beq
\sum_{j=-\infty}^\infty (-1)^j z^jq^{\bin{j+1}{2}} = (q)_\infty[qz;q]_\infty. 
\label{eq:1.6}
\eeq
As for the formulas \eqn{1.3} and \eqn{1.4}, they are, essentially, items G(4) and E(3), respectively, 
in Slater's table \cite{SL1}.

Several years ago, Andrews \cite{A1} revisited the umbral methods used by L.J. Rogers. 
In \cite{A1}, he discussed multidimensional identities of the form
\beqs
\sum_{\mathbf i}A(\mathbf L,q)\qBin{2\mathbf L}{\mathbf L+\mathbf i}{q}
=F(\mathbf L,q),
\eeqs
where $\mathbf L=(L_1,L_2,\ldots,L_d),\; \mathbf i=(i_1,i_2,\ldots,i_d)$ and
\beqs
\qBin{2\mathbf L}{\mathbf L+\mathbf i}{q}=\prod_{k=1}^d\qBin{2L_k}{L_k+i_k}{q}.
\eeqs
In particular, he proved that 
\beq
\begin{split}
\sum_{i,j,k}(-1)^{i+j+k} q^{\bin{i+j+k}{2}} & \qBin{2L}{L-i}{q} \qBin{2M}{M-j}{q} \qBin{2N}{N-k}{q} \\
= & \frac{(q)_{2L}(q)_{2M}(q)_{2N}}{(q)_{L+M-N}(q)_{L+N-M}(q)_{M+N-L}}, 
\end{split}
\label{eq:1.7}
\eeq
and that
\beq
\sum_{i,j}(-1)^i q^{\bin{i+j}{2}} \qBin{2L}{L-i}{q} \qBin{2M}{M-j}{q}
=(-1)^L\frac{(q;q^2)_L (q^2;q^2)_M (-1)_{M-L}}{(q)_{M-L}}.
\label{eq:1.8}
\eeq
We note that \eqn{1.7} is a three-dimensional generalization of \eqn{1.2}. 
Indeed, if we let $N=0$, we obtain that
\beq
\sum_{i,j}(-1)^{i+j} q^{\bin{i+j}{2}} \qBin{2L}{L-i}{q} \qBin{2M}{M-j}{q}=(q)_{2L}\delta_{L,M}.
\label{eq:1.9}
\eeq
If we now set $M=0$ we end up with \eqn{1.2}, as claimed. Also, \eqn{1.8} with $M=0$ reduces to \eqn{1.2}. 
On the other hand, \eqn{1.8} with $L=0$ becomes \eqn{1.5} with $z=-1,\;a=0$.

In what follows, we will use a small variant of \eqn{1.9}
\beq
\sum_{i,j}(-1)^{i+j} q^{\bin{i+j+1}{2}} \qBin{2L+1}{L-i}{q} \qBin{2M+1}{M-j}{q}=-(q)_{2L+1}\delta_{L,M}.
\label{eq:1.10}
\eeq
To verify that \eqn{1.10} holds for $M=0$ (or $L=0$) we use \eqn{1.5} twice as follows
\beqs
\begin{split}
& \sum_{i,j}(-1)^{i+j} q^{\bin{i+j+1}{2}}\qBin{2L+1}{L-i}{q} \qBin{1}{-j}{q} =\\
& \sum_i(-1)^i q^{\bin{i+1}{2}}\qBin{2L+1}{L-i}{q} - \sum_i(-1)^i q^{\bin{i}{2}}\qBin{2L+1}{L-i}{q}= \\
& (1)_{L+1} (q)_L-(q)_{L+1}(1)_L = 0-(q)_1\delta_{L,0}=-(q)_{2L+1}\delta_{L,0}. 
\end{split}
\eeqs
In \cite{R}, Riese used his {\it $q$MultiSum} package to provide a simple recurrence proof of \eqn{1.7}.
In the next section, I will rederive and generalize Riese's recurrences. 
As a bonus, I will get a uniform proof of \eqn{1.7} -- \eqn{1.10}. 
Moreover, I will show that the same proof can be employed to establish four new identities:
\beq
\begin{split}
\sum_{i,j} (-1)^i q^{(i+j)^2} & \qBin{2L}{L-i}{q^2} \qBin{2M}{M-j}{q^2} \\ 
& = (-1)^M \frac{(q;q^2)_{L-M}}{(-q;q^2)_{L-M}} (q^2;q^4)_L (q^2;q^4)_M, 
\end{split}
\label{eq:1.11}
\eeq

\beq
\sum_{i,j} (-1)^{i+j} q^{\bin{i+j}{2}} \qBin{2L}{L-i}{q} \qBin{2M}{M-j}{q^2} 
= (-1)^{L+M} q^{(L-M)^2} \frac{(q)_{2M}}{(q)_{M-L}}, 
\label{eq:1.12} 
\eeq

\beq
\begin{split}
\sum_{i,j} (-1)^{i+j} q^{2\bin{i+j}{2}} & \qBin{2L}{L-i}{q} \qBin{2M}{M-j}{q^2} \\
& = q^{L-M}(q;q^2)_L (-q;q^2)_{2M}(q^{2(L+1-M)};q^2)_M,
\end{split}
\label{eq:1.13} 
\eeq
and
\beq
\begin{split}
\sum_{i,j,k} (-1)^{i+j+k} q^{\bin{i+j+k}{2}} & \qBin{2L}{L-i}{q^2} \qBin{2M}{M-j}{q^2} \qBin{2N}{N-k}{q^2} \\
& =(q;q^2)_{L+M-N} (q;q^2)_{L+N-M}(q;q^2)_{M+N-L}.
\end{split}
\label{eq:1.14}
\eeq
We remark that \eqn{1.14} is a perfect quadratic analogue of \eqn{1.7}.
Setting $L=0$ in \eqn{1.11} yields \eqn{1.5} with $q\rightarrow q^2$ and $a=0,\;z=-\frac{1}{q}$.
Setting $M=0$ there yields \eqn{1.5} with $q\rightarrow q^2$, $a=0,\;z=\frac{1}{q}$.
Analogously, we can verify that \eqn{1.12} reduces to \eqn{1.2} and \eqn{1.3} and that \eqn{1.13} 
reduces to \eqn{1.4} and \eqn{1.5} with $q\rightarrow q^2,\;a=0,\;z=1$. 
Finally, \eqn{1.14} is a three-dimensional generalization of \eqn{1.3}.
Indeed, if we let $N=0$ in \eqn{1.14} we get
\beq
\sum_{i,j} (-1)^{i+j} q^{\bin{i+j}{2}} \qBin{2L}{L-i}{q^2} \qBin{2M}{M-j}{q^2} = (q;q^2)_{L+M} (-1)^{M+L} q^{(L-M)^2}.
\label{eq:1.15} 
\eeq
Letting $M=0$ in \eqn{1.15} yields \eqn{1.3}.

The remainder of this manuscript is organized as follows. 
In the next section, I will show how to use a simple polynomial identity to prove \eqn{1.7} -- \eqn{1.15}. 
In Section 3, I will employ some well-known $q$-binomial transformations to derive new two and three-dimensional identities. In Section 4, I will prove, among other results, that for $k\geq 0$
\beq
\begin{split}
\sum_{n_1,\ldots,n_{k+1}\geq 0}  & \frac{q^{N^2_1+N^2_2+\cdots+N^2_{k+1}}}{(q)_{n_1}(q)_{n_2}\cdots(q)_{n_k}
(q^2;q^2)_{n_{k+1}}(q;q^2)_{n_k+n_{k+1}}} \\
= & \frac{(q^{4+6k};q^{4+6k})_\infty}{(q)_\infty}\frac{[q^{2+2k};q^{4+6k}]_\infty}{[q^{1+k};q^{4+6k}]_\infty},
\end{split}
\label{eq:1.16} 
\eeq
where $N_i:=n_i+n_{i+1}+\cdots+n_{k+1},\;1\leq i\leq k+1$.\\
Finally, in Section 5, I will provide fresh insights into the Tri-Pentagonal
\begin{thm}{(Andrews)}
\beq
\frac{1}{(q)^3_\infty} \sum_{i,j,k} (-1)^{i+j+k}q^{\bin{i+j+k}{2}+i^2+j^2+k^2}
=\sum_{i,j,k\geq 0} \frac{q^{i^2+j^2+k^2}}{(q)_{i+j-k} (q)_{i+k-j}(q)_{j+k-i}}.
\label{eq:1.17} 
\eeq
\label{thm:1}
\end{thm}
Recently, this theorem was given a very interesting partition theoretical interpretation in \cite{A3}.

\bigskip
\section{\bf Simple polynomial identity and its implications}
\label{sec:2}
\medskip

It is a fair statement that one does not need a computer to check that 
\beq
\begin{split}
(1-x_1)(1-x_1q) & + (1-x_2)(1-x_2q) \\
& -(1-x_1x_2q)(1-x_1x_2q^2) \\
& - (1-x_1)(1-x_1q)(1-x_2)(1-x_2q) \\
& + x_1x_2(1+q)(1-x_1q)(1-x_2q)=0
\end{split}
\label{eq:2.1} 
\eeq
holds true for any $x_1,x_2,q \in \mathbb R$.\\
Next, we divide \eqn{2.1} by $(1-x_1x_2q)(1-x_1x_2q^2)$ and let $x_1=q^{L+a+i},\;x_2=q^{L-i}$ to get
\beq
\begin{split}
&  \frac{(1-q^{L+a+i})(1-q^{L+a+i+1})}{(1-q^{2L+a+1})(1-q^{2L+a+2})} \\
&+ \frac{(1-q^{L-i})(1-q^{L-i+1})}{(1-q^{2L+a+1})(1-q^{2L+a+2})} \\
&- \frac{(1-q^{L+a+i})(1-q^{L+a+i+1})(1-q^{L-i})(1-q^{L-i+1})}{(1-q^{2L+a+1})(1-q^{2L+a+2})} \\
&+ q^{2L+a}(1+q)\frac{(1-q^{L+a+i+1})(1-q^{L-i+1})}{(1-q^{2L+a+1})(1-q^{2L+a+2})}=1.
\end{split}
\label{eq:2.2} 
\eeq
If we multiply \eqn{2.2} by $\qBin{2(L+1)+a}{(L+1)-i}{q}$ we deduce that
\beq
\begin{split}
& \qBin{2L+a}{L-(i-1)}{q}+\qBin{2L+a}{L-(i+1)}{q}+q^{2L+a}(1+q)\qBin{2L+a}{L-i}{q} \\
& -(1-q^{2L+a})(1-q^{2L+a-1})\qBin{2(L-1)+a}{(L-1)-i}{q}=\qBin{2(L+1)+a}{(L+1)-i}{q}.
\end{split}
\label{eq:2.3} 
\eeq
We now define 
\beqs
f_a\bin{L,M}{i,j}=f_a\bin{L,M,q_1,q_2}{i,j}:=\qBin{2L+a}{L-i}{q_1}\qBin{2M+a}{M-j}{q_2}.
\eeqs
Obviously, \eqn{2.3} implies that 
\beq
\begin{split}
&   f_a\bin{L,M}{i-2,j}+f_a\bin{L,M}{i,j} + q_1^{2L+a}(1+q_1)f_a\bin{L,M}{i-1,j} \\
& -(1-q_1^{2L+a})(1-q_1^{2L+a-1})f_a\bin{L-1,M}{i-1,j}= f_a\bin{L+1,M}{i-1,j}
\end{split}
\label{eq:2.4} 
\eeq
and that
\beq
\begin{split}
&   f_a\bin{L,M}{i,j-2}+f_a\bin{L,M}{i,j} + q_2^{2M+a}(1+q_2)f_a\bin{L,M}{i,j-1} \\
& -(1-q_2^{2M+a})(1-q_2^{2M+a-1})f_a\bin{L,M-1}{i,j-1}= f_a\bin{L,M+1}{i,j-1}.
\end{split}
\label{eq:2.5} 
\eeq
It is easy to combine \eqn{2.4} and \eqn{2.5} as follows
\beq
\begin{split}
&   f_a\bin{L,M}{i-2,j} + q_1^{2L+a}(1+q_1)f_a\bin{L,M}{i-1,j} \\
& -(1-q_1^{2L+a})(1-q_1^{2L+a-1})f_a\bin{L-1,M}{i-1,j}-f_a\bin{L+1,M}{i-1,j} \\
& = f_a\bin{L,M}{i,j-2}+q_2^{2M+a}(1+q_2)f_a\bin{L,M}{i,j-1} \\
& -(1-q_2^{2M+a})(1-q_2^{2M+a-1})f_a\bin{L,M-1}{i,j-1}-f_a\bin{L,M+1}{i,j-1}.
\end{split}
\label{eq:2.6} 
\eeq
To proceed, we require one more definition
\beq
F_a(L,M)=F_a(L,M,x,y,q_0,q_1,q_2):=\sum_{i,j}x^iy^jq_0^{P(i+j)}f_a\bin{L,M,q_1,q_2}{i,j},
\label{eq:2.7} 
\eeq
where $P(z)$ is some polynomial in $z$.\\
Clearly,
\beq
y^2\sum_{i,j}x^iy^jq_0^{P(i+j-1)}f_a\bin{L,M}{i-2,j}=x^2\sum_{i,j}x^iy^jq_0^{P(i+j-1)}f_a\bin{L,M}{i,j-2}.
\label{eq:2.8} 
\eeq
Next, we multiply \eqn{2.6} by $x^iy^jq_0^{P(i+j-1)}$ and sum over $i,j$. 
Taking advantage of \eqn{2.8} we derive that
\beq
\begin{split}
&  x\bigl\{q_1^{2L+a}(1+q_1)F_a(L,M) -(1-q_1^{2L+a})(1-q_1^{2L+a-1})F_a(L-1,M) \\
& \qquad\qquad\qquad\qquad\qquad\quad - F_a(L+1,M)\bigr\} \\
& = y\bigl\{q_2^{2M+a}(1+q_2)F_a(L,M)-(1-q_2^{2M+a})(1-q_2^{2M+a-1})F_a(L,M-1) \\
& \qquad\qquad\qquad\qquad\qquad\quad - F_a(L,M+1)\bigr\},
\end{split}
\label{eq:2.9} 
\eeq
provided $|x|=|y|$.

Observe that for $a=0,1$ the recurrence \eqn{2.9} together with the boundary values $F_a(L,0),\;L\ge 0$ and
$F_a(0,M),\;M\ge 0$ specifies $F_a(L,M)$ completely for $L\ge 0,\;M\ge 0$.

It is plain that the left hand sides in \eqn{1.8}--\eqn{1.13} and \eqn{1.15} are of the form \eqn{2.7} with
$|x|=|y|=1$. It is also straightforward to check that the right hand sides there satisfy \eqn{2.9}. This 
implies that \eqn{1.8}--\eqn{1.13} and \eqn{1.15} hold true if they hold when $L\ge 0,\;M=0$ and $L=0,\;M\ge 0$.
But this is indeed the case as we saw in the Introduction.

Fortunately, more is true. Define
\beq
Y_m(L,M,N):=\sum_{i,j,k}(-1)^{i+j+k} q^{\bin{i+j+k}{2}} 
\qBin{2L}{L+i}{q^m} \qBin{2M}{M+j}{q^m} \qBin{2N}{N+k}{q^m}.
\label{eq:2.10} 
\eeq
Clearly, $Y_1$ and $Y_2$ are the left hand sides of \eqn{1.7} and \eqn{1.14}, respectively. 
Multiply \eqn{2.6} with $q_1=q_2=q^m$ by $(-1)^{i+j+k}q^{\bin{i+j+k-1}{2}}\qBin{2N}{N-k}{q^m}$ and
sum over $i,j,k$ to derive that
\beq
\begin{split}
&  q^{2Lm}(1+q^m)Y_m(L,M,N)-(1-q^{2Lm})(1-q^{2Lm-m})Y_m(L-1,M,N) \\
& \qquad\qquad\qquad\qquad\qquad\quad -Y_m(L+1,M,N) \\
& = q^{2Mm}(1+q^m)Y_m(L,M,N)-(1-q^{2Mm})(1-q^{2Mm-m})Y_m(L,M-1,N) \\
& \qquad\qquad\qquad\qquad\qquad\quad -Y_m(L,M+1,N). 
\end{split}
\label{eq:2.11} 
\eeq
We remark that \eqn{2.11} with $m=1$ is, essentially, the recurrence derived by Riese in \cite{R}.
Once again, \eqn{2.11} together with the boundary values $Y_m(0,M,N),\;M\ge 0,\;N\ge 0$ and 
$Y_m(L,0,N),\;L\ge 0,\;N\ge 0$ specify $Y_m(L,M,N)$ completely for $L,M,N\ge 0$. Next, we check that
the right hand sides of \eqn{1.7} and \eqn{1.14} satisfy \eqn{2.11} with $m=1$ and $m=2$, respectively.
Moreover, on the boundary these identities reduce to the two proven identities: \eqn{1.9} and \eqn{1.15}.
And so, \eqn{1.7} and \eqn{1.14} hold true, as claimed.

The reader may wonder if the polynomial identity \eqn{2.1} can be extended to $n$ variables 
$:x_1,x_2,x_3,\ldots,x_n$. This is indeed possible. The following generalization was suggested to me by
Alain Lascoux:
\beq
\prod_{i=1}^n(x_i)_2+\sum_{t=1}^n\sum_{i_1<i_2<\cdots<i_t}(-1)^t(q^{t-1}x_{i_1}x_{i_2}\cdots x_{i_t})_2 
= \prod_{i=1}^n(qx_i)_1\sum_{t=1}^{n-1}(q^t-(-1)^t)e_{t+1},
\label{eq:2.12} 
\eeq
where $e_i$'s are the elementary symmetric functions in $x_1,x_2,\ldots,x_n$.

\bigskip
\section{\bf $q$-binomial transformations}
\label{sec:3}
\medskip

We begin by recalling some well-known formulas
\beq
\sum_{r\ge 0}q^{r^2}\frac{(q)_{2L}}{(q)_{L-r}(q)_{2r}}\qBin{2r}{r-j}{q}=q^{j^2}\qBin{2L}{L-j}{q},
\label{eq:3.1} 
\eeq
\beq
\sum_{r\ge 0}q^{L-2r}(-\frac{1}{q};q^2)_{L-2r} \qBin{L}{2r}{q^2} \qBin{2r}{r-j}{q^4}=\qBin{2L}{L-2j}{q},
\label{eq:3.2} 
\eeq
\beq
\sum_{r\ge 0}q^{2r(r+a)}(-q)_{L-2r-a} \qBin{L}{2r+a}{q^2} \qBin{2r+a}{r-j}{q^2}=q^{2j(j+a)}\qBin{2L}{L-2j-a}{q},
\label{eq:3.3} 
\eeq
where $a=0,1$. \\
We remark that \eqn{3.1} was used by Bressoud \cite{Br} to give a simple proof of the Rogers--Ramanujan identities.
It can be recognized as a special case of the Bailey Lemma in its version due to Andrews \cite{A2} and Paule \cite{P}.
The transformations \eqn{3.2} and \eqn{3.3} were introduced by Berkovich and Warnaar in \cite{BW}.

In \cite{A1}, Andrews applied \eqn{3.1} to \eqn{1.7} three times to obtain the tri-pentagonal theorem \eqn{1.17}.
It is interesting that a single application of \eqn{3.1} to \eqn{1.7} yields a new three-dimensional identity 
\beq
\begin{split}
& \sum_{i,j,k} (-1)^{i+j+k} q^{\bin{i+j+k}{2}+i^2} \qBin{2L}{L-i}{q} \qBin{2M}{M-j}{q} \qBin{2N}{N-k}{q} \\
& =q^{(N-M)^2}(q)_{L+M+N}\qBin{2L}{L+N-M}{q}.
\end{split}
\label{eq:3.4}
\eeq
Indeed, we have that
\beq
\begin{split}
& \sum_{i,j,k} (-1)^{i+j+k} q^{\bin{i+j+k}{2}+i^2} \qBin{2L}{L-i}{q} \qBin{2M}{M-j}{q} \qBin{2N}{N-k}{q} \\
& = \sum_{r\ge 0}\frac{q^{r^2}}{(q)_{L-r}}\frac{(q)_{2L}(q)_{2M}(q)_{2N}}{(q)_{r+M-N}(q)_{r+N-M}(q)_{M+N-r}}.
\end{split}
\label{eq:3.5}
\eeq
The right hand side of \eqn{3.5} can be written in the form
\beq
\begin{split}
\mbox{RHS\eqn{3.5}} = \frac{q^{n^2}(q)_{2L}(q)_{2M}(q)_{2N}}{(q)_{2n}(q)_{L-n}(q)_{M+N-n}}\theta(L\ge n)
\sum_{r\ge 0}\frac{(q^{n-L},q^{n-M-N})_r}{(q,q^{2n+1})_r} q^{r(L+M+N+1)},
\end{split}
\label{eq:3.6}
\eeq
where $n:=|N-M|$ and 
\begin{align}
\theta(L\ge n)=
\begin{cases} 1, & \mbox{if } L\ge n, \\
              0, & \mbox{otherwise}. 
\end{cases}
\label{eq:3.7}
\end{align}
The sum in \eqn{3.6} can be evaluated by the $q$-Chu-Vandermonde formula [\cite{GR}, (II.7)] as
\beqs
\frac{(q)_{M+N+L}}{(q^{2n+1})_{M+N-n}(q)_{n+L}}.
\eeqs
And so,
\beqs
\begin{split}
RHS\eqn{3.5} & =\frac{q^{n^2}(q)_{2L}(q)_{2M}(q)_{2N}(q)_{M+N+L}}{(q)_{M+N+n}(q)_{M+N-n}(q)_{L-n}(q)_{L+n}} \\
& =q^{(N-M)^2}\qBin{2L}{L+N-M}{q}(q)_{L+M+N},
\end{split}
\eeqs
as claimed.

Obviously, \eqn{3.2} and \eqn{3.3} can also be employed to produce new two-dimensional identities. 
For example, we can replace $q$ by $q^4$ in \eqn{1.9} and apply \eqn{3.2} to obtain
\beq
\sum_{i,j} (-1)^{i+j} q^{4\bin{i+j}{2}} \qBin{2L}{L-i}{q^4} \qBin{2M}{M-2j}{q}
= q^{M-2L}(-\frac{1}{q};q^2)_{M-2L} \qBin{M}{2L}{q^2} (q^4;q^4)_{2L}.  
\label{eq:3.8}
\eeq
Or we can replace $q$ by $q^2$ in \eqn{1.9} and apply \eqn{3.3} with $a=0$ twice to get
\beq
\begin{split}
& \sum_{i,j} (-1)^{i+j} q^{(i+j)(i+j+1)+2i^2+2j^2} \qBin{2L}{L-2i}{q} \qBin{2M}{M-2j}{q} \\
&=\sum_{\substack{r\equiv 0\mod 2,\\r\ge0}}q^{r^2}\frac{(q^2;q^2)_L (q^2;q^2)_M}{(q^2;q^2)_r(q)_{L-r}(q)_{M-r}}.
\end{split}
\label{eq:3.9}
\eeq
Analogously, replacing $q$ by $q^2$ in \eqn{1.10} and using \eqn{3.3} with $a=1$ we obtain that
\beq
\begin{split}
& -q\sum_{i,j} (-1)^{i+j} q^{(i+j)(i+j+1)+2i(i+1)+2j(j+1)} \qBin{2L}{L-(2i+1)}{q} \qBin{2M}{M-(2j+1)}{q} \\
&=\sum_{\substack{r\equiv 1\mod 2,\\r>0}} q^{r^2} \frac{(q^2;q^2)_L(q^2;q^2)_M}{(q^2;q^2)_r(q)_{L-r}(q)_{M-r}}.
\end{split}
\label{eq:3.10}
\eeq

\bigskip
\section{\bf Two infinite families of multiple series identities}
\label{sec:4}
\medskip

If we let $L\rightarrow\infty$ in \eqn{3.9} and \eqn{3.10} we end up with the following result
\beq
\begin{split}
& \sum_{\substack{r\equiv a\mod 2,\\r\ge0}} q^{r^2} \frac{(q^2;q^2)_M}{(q^2;q^2)_r(q)_{M-r}}=\frac{(-q)^a}{(q^2;q^2)_\infty} \\
& \sum_{i,j} (-1)^{i+j} q^{(i+j)(i+j+1)+2i(i+a)+2j(j+a)} \qBin{2M}{M-2j-a}{q},
\end{split}
\label{eq:4.1}
\eeq
where $a=0,1$.
Remarkably, the double sum on the right hand side of \eqn{4.1} can be reduced to a single sum. 
To this end, we perform a clever change of the summation variables $j\rightarrow 3j+r-a$ with $r=0,\pm 1$ and
$i\rightarrow i-j$. This yields 
\beq
\begin{split}
& RHS\eqn{4.1}= \frac{1}{(q^2;q^2)_\infty}\sum_{r=-1}^1 (-1)^r\sum_{i=-\infty}^\infty (-1)^i q^{3i^2+i(1+2r)} \\
& \sum_{j=-\infty}^\infty q^{24j^2+2j(8r+1-6\delta_{a,1})+3r^2+r+\delta_{a,1}(1-4r)} \qBin{2M}{M-6j-2r+a}{q}.
\end{split}
\label{eq:4.2}
\eeq
We now make use of a special case of \eqn{1.6} 
\beq
\sum_{i=-\infty}^\infty (-1)^i q^{3i^2+i(1+2r)}=(q^2;q^2)_\infty(1-\delta_{r,1}),\quad r=0,\pm1
\label{eq:4.3}
\eeq
to simplify \eqn{4.2} further. This way we obtain
\beq
\begin{split}
& \sum_{\substack{r\equiv a\mod 2,\\r\ge0}} q^{r^2} \frac{(q^2;q^2)_M}{(q^2;q^2)_r(q)_{M-r}}= \\
& \sum_{j=-\infty}^\infty q^{24j^2+2j(1-6\delta_{a,1})+\delta_{a,1}} \qBin{2M}{M-6j+a}{q} - \\
& \sum_{j=-\infty}^\infty q^{24j^2-2j(7+6\delta_{a,1})+5\delta_{a,1}} \qBin{2M}{M-6j+a+2}{q}.
\end{split}
\label{eq:4.4}
\eeq
Clearly, one could have arrived at \eqn{4.4} by taking a more direct route by applying \eqn{3.3} to the 
polynomial identity
\beqs
\sum_{j=-\infty}^\infty q^{(3j+1+a)2j}\left(\qBin{2M+a}{M-3j}{q^2} - \qBin{2M+a}{M-3j-1}{q^2}\right)=q^{2M(M+a)},
\quad a=0,1,
\eeqs
which is, essentially, A(5) and A(8) in \cite{SL1}. \\
However, I feel that the passage from \eqn{4.1} to \eqn{4.4} is a good warm-up exercise to prepare the reader for the development in the next section.

We can now follow a well trodden path \cite{A2} and iterate \eqn{3.1} to get for $k\ge 1$ and $a=0,1$
\beq
\begin{split}
& \sum_{\substack{n_1,\ldots,n_{k+1}\ge0,\\n_{k+1}\equiv a\mod 2}}
\frac{q^{N^2_1+\cdots+N^2_{k+1}}(q)_{2M}}{(q)_{M-N_1}(q)_{n_1}\cdots(q)_{n_k}(q^2;q^2)_{n_{k+1}}(q;q^2)_{n_k+n_{k+1}}}= \\
& \sum_{j=-\infty}^\infty q^{12(2+3k)j^2+2j(1-\delta_{a,1}6(1+k))+(1+k)\delta_{a,1}} \qBin{2M}{M-6j+a}{q} - \\
& \sum_{j=-\infty}^\infty q^{12(2+3k)j^2+2j(-7-12k-\delta_{a,1}6(1+k))+2+k(a+2)^2+5\delta_{a,1}}
\qBin{2M}{M-6j+a+2}{q},
\end{split}
\label{eq:4.5}
\eeq
where $N_i=n_i+\cdots+n_{k+1},\;1\le i\le k+1$.
If we now let $M\rightarrow\infty$ in \eqn{4.5} we arrive at
\beq
\begin{split}
& \sum_{\substack{n_1,\ldots,n_{k+1}\ge0,\\n_{k+1}\equiv a\mod 2}}
\frac{q^{N^2_1+\cdots+N^2_{k+1}}}{(q)_{n_1}\cdots(q)_{n_k}(q^2;q^2)_{n_{k+1}}(q;q^2)_{n_k+n_{k+1}}}= \\
& \frac{q^{(1+k)\delta_{a,1}}}{(q)_\infty} \sum_{j=-\infty}^\infty q^{4(2+3k)(3j-1)j} z_a^{3j}(1-z_aq^{8(2+3k)j}),
\end{split}
\label{eq:4.6}
\eeq
where $a=0,1$ and $z_a:=q^{2+4a+4k(1+\delta_{a,1})}$.
At this stage we recall the quintuple product identity [\cite{GR}, Ex. 5.6]
\beq
\begin{split}
\sum_{n=-\infty}^\infty (-1)^n q^{\frac{3n-1}{2}n}z^{3n}(1+zq^n)& =(q,-z,-\frac{q}{z})_\infty (qz^2,\frac{q}{z^2};q^2)_\infty \\ 
& = (q)_\infty \frac{[z^2;q]_\infty}{[z;q]_\infty}. 
\end{split}
\label{eq:4.7}
\eeq
This identity enables us to rewrite \eqn{4.5} as 
\beq
\begin{split}
& \sum_{\substack{n_1,\ldots,n_{k+1}\ge0,\\n_{k+1}\equiv a\mod 2}}
\frac{q^{N^2_1+\cdots+N^2_{k+1}}}{(q)_{n_1}\cdots(q)_{n_k}(q^2;q^2)_{n_{k+1}}(q;q^2)_{n_k+n_{k+1}}}= \\
& q^{(1+k)\delta_{a,1}}\frac{(q^{16+24k};q^{16+24k})_\infty}{(q)_\infty}
[z_a;q^{16+24k}]_\infty [q^{16+24k}z_a^2;q^{32+48k}]_\infty,
\end{split}
\label{eq:4.8}
\eeq
where $a=0,1$.

It remains to establish \eqn{1.16}. To this end we add together \eqn{4.5} with $a=0$ and \eqn{4.5} with $a=1$. 
This way we immediately obtain the correct left hand side of \eqn{1.16}. 
Making use of \eqn{4.7} on the right we derive that
\beq
\begin{split}
& \frac{1}{(q)_\infty} \sum_{j=-\infty}^\infty q^{4(2+3k)(3j-1)j} z_1^{3j}(1-z_1q^{8(2+3k)j}) + \\
& \frac{q^{1+k}}{(q)_\infty} \sum_{j=-\infty}^\infty q^{4(2+3k)(3j-1)j} z_2^{3j}(1-z_2q^{8(2+3k)j}) = \\
& \frac{1}{(q)_\infty}
\sum_{j=-\infty}^\infty (-1)^j q^{(2+3k)(3j-1)j} q^{3(1+k)j} (1+q^{1+k}q^{(2+3k)2j})= \\
& \qquad\qquad\qquad \frac{(q^{4+6k};q^{4+6k})_\infty}{(q)_\infty} \frac{[q^{2+2k};q^{4+6k}]_\infty}{[q^{1+k};q^{4+6k}]_\infty},
\end{split}
\label{eq:4.9}
\eeq
as desired. It is instructive to compare \eqn{1.16} with a somewhat similar formula [\cite{A2}, (1.8)]
\beq
\begin{split}
\sum_{n_1,\ldots,n_k\ge 0}& \frac{q^{\tilde N^2_1+\cdots+\tilde N^2_k}}
{(q)_{n_1}(q)_{n_2}\cdots(q)_{n_{k-1}}(q)_{2n_k}}=\\
& \frac{(q^{4+6k};q^{4+6k})_\infty}{(q)_\infty} \frac{[q^{2+2k};q^{4+6k}]_\infty}{[-q^{1+k};q^{4+6k}]_\infty},
\end{split}
\label{eq:4.10}
\eeq
where $\tilde N_i:=n_i+\cdots+n_k,\;1\le i\le k$.

\bigskip
\section{\bf One-dimensional version of the tri-pentagonal theorem}
\label{sec:5}
\medskip

This paper arose from my attempt to ascertain if Andrews' formula \eqn{1.17} was ``genuinely" three-dimensional. 
In the course of this investigation I found that \eqn{1.17} can be ``flattened" as follows
\beq
\begin{split}
\frac{(q^4;q^4)_\infty(q^{20};q^{20})_\infty}{(q)^2_\infty}
& \sum^1_{a=0} (-q)^a [-q^{2-2a};q^4]_\infty [q^{4+2a};q^{20}]_\infty [q^{12-4a};q^{40}]_\infty = \\
\frac{1}{2}(-\sqrt q)_\infty 
& \sum_{i\ge 0}\frac{(-\sqrt q)_i q^{\frac{3}{2}i^2}}{(q)_{2i}} + 
\frac{1}{2}(\sqrt q)_\infty\sum_{i\ge 0}\frac{(\sqrt q)_i q^{\frac{3}{2}i^2}(-1)^i}{(q)_{2i}}.
\end{split}
\label{eq:5.1}
\eeq
To reduce the triple sum on the left of \eqn{1.17} to a single sum on the left of \eqn{5.1}, we perform some clever changes $j\rightarrow 3j-i+r$ with $r=0,\pm1$ and $k\rightarrow k-j$. Using \eqn{4.3} with $q^2\rightarrow q$, 
we obtain that 
\beqs
LHS\eqn{1.17}= \frac{1}{(q)^2_\infty} \sum_{i,j} q^{12j^2+j+2i(i-3j)} 
- \frac{1}{(q)^2_\infty} \sum_{i,j} q^{12j^2-7j+2i(i-3j)+2i+1}.
\eeqs
Next, we cleverly substitute $2j-a$ for $j$ and $i+3j-a$ for $i$  with $a=0,1$ in the first sum and with $a=0,-1$ 
in the second sum. This way we get
\beqs
\begin{split}
LHS\eqn{1.17} & = \frac{1}{(q)^2_\infty} \sum_{i,j} \left\{q^{30j^2+2j+2i^2}-q^{30j^2-8j+2i^2+2i+1}\right. \\
& + \left. q^{30j^2-28j+2i^2+2i+7}-q^{30j^2+22j+2i^2+4}\right\}\\
& = \frac{1}{(q)^2_\infty}\sum_{i=-\infty}^\infty q^{2i^2}\sum_{j=-\infty}^\infty q^{30j^2+2j}(1-q^4q^{20j}) \\
& - \frac{q}{(q)^2_\infty}\sum_{i=-\infty}^\infty q^{2i^2+2i}\sum_{j=-\infty}^\infty q^{30j^2+8j}(1-q^6q^{20j}) \\
& = LHS\eqn{5.1},
\end{split}
\eeqs
where we used \eqn{1.6} and \eqn{4.7} in the last step.
Next, \eqn{3.4} with $L\rightarrow\infty$ suggests that
\beq
RHS\eqn{1.17}= \sum_{i,j\ge0} \frac{q^{i^2+j^2+(i-j)^2}}{(q)_{2i}(q)_{2j}}. 
\label{eq:5.2}
\eeq
In other words, one variable on the right of \eqn{1.17} can be summed out by the $q$-Chu-Vandermonde formula as follows.
First, replace $k$ by $k+n$ on the right of \eqn{1.17}, where $n=|i-j|$. This leads to
\beqs
\begin{split}
RHS\eqn{1.17} & = \sum_{i,j,k\ge0} \frac{q^{i^2+j^2+(k+n)^2}}{(q)_k(q)_{k+2n}(q)_{i+j-n-k}}\\
& = \sum_{i,j\ge0} \frac{q^{i^2+j^2+n^2}}{(q)_{2n}(q)_{i+j-n}} \sum_{k\ge0} \frac{(q^{n-i-j})_k}{(q,q^{1+2n})_k}
q^{\bin{k}{2}}q^{k(1+i+j+n)}(-1)^k.
\end{split}
\eeqs
The inner sum on the right can be evaluated by [\cite{GR}, (II.7)] with $a\rightarrow\infty$ as 
\beqs
\frac{1}{(q^{1+2n})_{i+j-n}}. 
\eeqs
Finally, using
\beqs
(q)_{i+j-n}(q)_{2n}(q^{1+2n})_{i+j-n}=(q)_{i+j-n}(q)_{i+j+n}=(q)_{2i}(q)_{2j},
\eeqs
we end up with \eqn{5.2}.
Actually, one can sum out $j$ in \eqn{5.2}, as well. This can be done as follows.
\beqs
\begin{split}
RHS\eqn{1.17}& = \sum_{i\ge0}\frac{q^{2i^2}}{(q)_{2i}} 
\sum_{j\ge0}\frac{q^{\bin{j}{2}}}{(q)_j}q^{(\frac{1}{2}-i)j}\frac{1+(-1)^j}{2} \\
& = \frac{1}{2}\sum_{i\ge0}\frac{q^{2i^2}}{(q)_{2i}}
\left((-q^{\frac{1}{2}-i})_\infty +(q^{\frac{1}{2}-i})_\infty\right), 
\end{split}
\eeqs
where we have made use of the Euler identity [\cite{GR}, (II.2)]
\beqs
\sum_{j\ge0}\frac{q^{\bin{j}{2}}}{(q)_j}z^j=(-z)_\infty. 
\eeqs
It is not hard to verify that
\beqs
(\pm q^{\frac{1}{2}-i})_\infty=(\mp1)^i q^{-\frac{i^2}{2}}(\pm\sqrt q)_i (\pm\sqrt q)_\infty. 
\eeqs
And so,
\beqs
RHS\eqn{1.17}=RHS\eqn{5.1},
\eeqs
as claimed. Thus, we have shown that \eqn{1.17} and \eqn{5.1} are, indeed, equivalent.

I would like to finish this section by providing an independent proof of \eqn{5.1}. This proof is based on the following identities
\beq
\sum_{i\ge 0}\frac{q^{3i^2}(q;q^2)_i(-1)^i}{(q^2;q^2)_{2i}}= \frac{(q;q^2)_\infty(q^4;q^4)_\infty}{(q^2;q^2)_\infty}
\sum_{n\ge 0}\frac{q^{n^2}}{(q^4;q^4)_n}
\label{eq:5.3}
\eeq
and
\beq
\sum_{i\ge 0}\frac{q^{i^2+ai}}{(q)_{2i+a}}= \frac{(q^{10};q^{10})_\infty}{(q)_\infty}
[q^{2+a};q^{10}]_\infty [q^{6-2a};q^{20}]_\infty,
\label{eq:5.4}
\eeq
where $a=0,1$. These identities are in Slater's list: \eqn{5.3} is a combination of items (19) and (20) in \cite{SL2} and \eqn{5.4} with $a=0,1$ are items (98) and (94), respectively. Given a power series $f(q)$, let
\beq
\{f(q)\}_e:=\frac{f(q)+f(-q)}{2}.
\label{eq:5.5}
\eeq
Substituting $q^2$ for $q$ in \eqn{5.1} and making use of \eqn{5.3} and \eqn{5.5}, we can rewrite \eqn{5.1} in the form 
\beq
\begin{split}
& \left\{ \frac{(q^2,q,q;q^2)_\infty}{(q^2;q^2)^2_\infty}(q^4;q^4)_\infty
\sum_{n\ge 0}\frac{q^{n^2}}{(q^4;q^4)_n}\right \}_e = \\
& \frac{(q^8;q^8)_\infty (q^{40};q^{40})_\infty}{(q^2;q^2)^2_\infty} 
\sum^1_{a=0}(-q^2)^a [-q^{4-4a};q^8]_\infty [q^{8+4a};q^{40}]_\infty [q^{24-8a};q^{80}]_\infty.
\end{split}
\label{eq:5.6}
\eeq
Finally, we use \eqn{1.6} on the left of \eqn{5.6} to deduce that
\beqs
\begin{split}
& LHS\eqn{5.6} =\\
&\frac{(q^4;q^4)_\infty}{(q^2;q^2)^2_\infty}
\left\{ \sum^\infty_{j=-\infty}(-1)^j q^{j^2} \sum_{n\ge 0}\frac{q^{n^2}}{(q^4;q^4)_n} \right\}_e = \\
& \frac{(q^4;q^4)_\infty}{(q^2;q^2)^2_\infty} 
\left\{ \sum^\infty_{j=-\infty}q^{4j^2} \sum_{n\ge 0}\frac{q^{4n^2}}{(q^4;q^4)_{2n}}
- q^2\sum^\infty_{j=-\infty}q^{4j^2+4j} \sum_{n\ge 0}\frac{q^{4n^2+4n}}{(q^4;q^4)_{2n+1}}\right\} \\ 
& = RHS\eqn{5.6},
\end{split}
\eeqs
where we employed \eqn{1.6} and \eqn{5.4} with $q\rightarrow q^4$ in the last step.\\
And so, \eqn{5.6} and, consequently, \eqn{5.1} hold true.

\bigskip
\subsection*{Acknowledgment}
I would like to thank George Andrews, who for more than two years kept insisting that I write up these humble observations. I am grateful to Axel Riese for many constructive discussions and to Michael Somos for his careful reading of this manuscript.


\end{document}